\documentclass{article}
\usepackage[round]{natbib}
\usepackage{mathrsfs}
\usepackage{latexsym}
\usepackage{amsmath}
\usepackage{amsfonts}
\usepackage{url}
\usepackage{amssymb}
\usepackage{amsthm}
\usepackage{graphicx} %
\usepackage{subfigure} 
\usepackage{epstopdf}
\usepackage{float}
\usepackage{natbib}
\newcommand{\inv}{^{-1}} %
\newcommand{\TV}[1]{\nrm{#1}_{\textrm{{\tiny \textup{TV}}}}}

\usepackage{dsfont}
\newcommand{\chr}{\mathds{1}}
\newcommand{\pred}[1]{\chr_{\left\{ #1 \right\}}}

\renewcommand{\vec}[1]{\ensuremath{\text{{\bf\textrm{#1}}}}}
\newcommand{\E}{\mathbb{E}}
\newcommand{\V}{\operatorname{Var}}

\newcommand{\sgn}{\operatorname{sign}}

\renewcommand{\P}{\mathbb{P}}

\newcommand{\ben}{\begin{enumerate}}
\newcommand{\een}{\end{enumerate}}
\newcommand{\bit}{\begin{itemize}}
\newcommand{\eit}{\end{itemize}}

\renewcommand{\vec}[1]{\bs{\mathrm{#1}}}

\def\clap#1{\hbox to 0pt{\hss#1\hss}}

\newcommand{\nrm}[1]{\left\Vert #1 \right\Vert}

\newcommand{\R}{\mathbb{R}}

\newcommand{\beq}{\begin{eqnarray*}}
\newcommand{\eeq}{\end{eqnarray*}}
\newcommand{\beqn}{\begin{eqnarray}}
\newcommand{\eeqn}{\end{eqnarray}}
\newcommand{\paren}[1]{\left( #1 \right)}
\newcommand{\sqprn}[1]{\left[ #1 \right]}
\newcommand{\tlprn}[1]{\left\{ #1 \right\}}
\newcommand{\set}[1]{\tlprn{#1}}
\newcommand{\abs}[1]{\left| #1 \right|}

\newcommand{\gn}{\, | \,}
\newcommand{\ds}{\displaystyle}
\newcommand{\ts}{\textstyle}
\newcommand{\bs}{\boldsymbol}
\renewcommand{\th}{\ensuremath{^{\mathrm{th}}}~}
\def\longto{\mathop{\longrightarrow}\limits}

\newcommand{\hide}[1]{}
\newcommand{\oo}[1]{\frac{1}{#1}}
\def\eps{\varepsilon}

\newtheorem{theorem}{Theorem}
\newtheorem{corollary}[theorem]{Corollary}

\newtheorem{lemma}[theorem]{Lemma}
\newtheorem{remark}{Remark}
\newtheorem*{remark*}{Remark}
\newcommand{\bepf}{\begin{proof}}
\newcommand{\enpf}{\end{proof}}

\newcommand{\hideshow}[2]{#2}
\newcommand{\shortlong}[2]{{#2}}
\newcommand{\CC}{} %
\title{
Consistency of weighted majority votes
}
\author{
Daniel Berend \and
Aryeh Kontorovich
}
\begin{document}
\maketitle
\begin{abstract}
We revisit the classical decision-theoretic problem of weighted expert voting
from a statistical learning perspective.
In particular, we examine the consistency (both asymptotic and finitary)
of the optimal Nitzan-Paroush weighted majority
and related rules. In the case of known
expert competence levels, we give 
sharp error estimates for the optimal rule.
When the competence levels are unknown,
they must be empirically estimated.
We provide frequentist and Bayesian analyses for this situation.
Some of our proof techniques are non-standard and may be of 
independent interest. 
The bounds we derive are nearly optimal,
and several challenging open problems are posed.
Experimental results are provided to illustrate the theory.
\end{abstract}
\newcommand{\weta}{\vec w\cdot\vec\eta}
\newcommand{\veta}{\vec\eta}
\newcommand{\hf}{{\ts\oo2}}
\newcommand{\implied}{\Longleftarrow}
\newcommand{\argmin}{\mathop{\mathrm{argmin}}}
\newcommand{\argmax}{\mathop{\mathrm{argmax}}}
\newcommand{\lc}{^{\textrm{{\tiny \textup{LC}}}}}
\newcommand{\hc}{^{\textrm{{\tiny \textup{HC}}}}}
\newcommand{\ba}{^{\textrm{{\tiny \textup{Ba}}}}}
\newcommand{\opt}{^{\textrm{{\tiny \textup{OPT}}}}}
\newcommand{\maj}{^{\textrm{{\tiny \textup{MAJ}}}}}
\newcommand{\tru}{^{\textrm{{\tiny \textup{TRU}}}}}
\newcommand{\Bdist}{\operatorname{Beta}}
\newcommand{\conv}{\operatorname{conv}}
\newcommand{\minpq}{\min\set{p_i,q_i}}
\newcommand{\LAND}{\cap}
\newcommand{\BLAND}{\bigcap}
\newcommand{\LOR}{\cup}
\newcommand{\BLOR}{\bigcup}
\section{Introduction}
The problem of weighting the input of several experts arises 
in many situations
and is of considerable theoretical and practical importance. 
The rigorous study of majority vote
has its roots in
the work of Condorcet \citeyearpar{de1785essai}.
By the 70s,
the field of decision theory was actively exploring 
various
voting rules
(see \citet{nitzan-paroush82} and the references therein).
A typical setting is as follows.
An agent is
tasked with predicting some random variable $Y\in\set{\pm1}$
based on 
input $X_i\in\set{\pm1}$ from each of $n$ experts.
Each expert $X_i$ 
has a {\em competence} level $p_i\in(0,1)$, 
which is the probability of making a correct
prediction:
$\P(X_i=Y)=p_i$.
Two simplifying assumptions are commonly made:
\begin{itemize}
\item[(i)] {\em Independence}: The random variables $\set{X_i:i\in[n]}$ are mutually independent.
\item[(ii)] {\em Unbiased truth}: $\P(Y=+1)=\P(Y=-1)=1/2$.
\end{itemize}
We will discuss these assumptions below in greater detail; for now, let us just take them as given.
(Since the bias of $Y$ can be easily estimated from data, only the independence assumption is truly restrictive.)
A {\em decision rule} is a mapping $f:\set{\pm1}^n\to\set{\pm1}$ 
from
the $n$ expert inputs to
the agent's final decision. 
Our 
quantity of interest
throughout the paper will be 
the agent's probability of error,
\beqn
\label{eq:prerr}
\P(f(\vec X)\neq Y).
\eeqn
\hide{}\hideshow{}{A decision rule $f$ is {\em optimal} if it minimizes
the quantity in (\ref{eq:prerr})
over all 
possible
decision rules.}
\citet{nitzan-paroush82} 
showed that, 
when Assumptions (i)--(ii) hold and the true competences $p_i$ are known,
the 
optimal decision rule
is obtained by an appropriately weighted majority vote:
\beqn
\label{eq:N-P}
f\opt(\vec x)=\sgn\paren{\sum_{i=1}^n
w_i
x_i},
\eeqn
where the weights $w_i$ are given by
\beqn
\label{eq:wdef}
w_i = \log\frac{p_i}{1-p_i},
\qquad i\in[n]
.
\eeqn
Thus, $w_i$ is the log-odds of expert $i$ being correct 
---
and the voting rule in (\ref{eq:N-P}),
also known as {\em naive Bayes} \citep{hastie2009the},
 may be seen as a simple consequence
of the Neyman-Pearson lemma \citep{NP1933}.
\paragraph{Main results.}
The formula in (\ref{eq:N-P}) raises immediate questions,
which apparently have not previously been addressed.
The first one has to do with the {\em consistency} of
the Nitzan-Paroush optimal rule: under what conditions does
the probability of error decay to zero
and 
at what rate?
In 
Section~\ref{sec:known},
we show that 
the probability of error is
controlled by the 
{\em committee potential} $\Phi$, defined by
\beqn
\label{eq:Phi}
\Phi=\sum_{i=1}^n (p_i-\hf)w_i
=\sum_{i=1}^n (p_i-\hf)\log\frac{p_i}{1-p_i}.
\eeqn
More precisely,
we prove in Theorem~\ref{thm:N-P-known} that
\beq
-\log \P(f\opt(\vec X)\neq Y)
&\asymp &
\Phi,
\eeq
where $\asymp$ denotes equivalence up to universal multiplicative constants.
Another issue not addressed by the Nitzan-Paroush result is how to handle the case
where the competences $p_i$ are not known exactly but
rather estimated empirically by $\hat p_i$.
We present two solutions to this problem: a frequentist and a Bayesian one.
As we show in
Section~\ref{sec:freq},
the frequentist approach does not admit an optimal empirical decision rule.
Instead, we analyze empirical decision rules in various settings:
high-confidence (i.e., $|\hat p_i-p_i|\ll1$) vs. low-confidence,
adaptive vs. nonadaptive.
The low-confidence regime requires no additional assumptions,
but provides weaker guarantees
(Theorem~\ref{thm:unknown-freq-low}).
In the high-confidence regime,
the adaptive approach provides error estimates in terms of
the empirical $\hat p_i$s (Theorem~\ref{thm:adapt}),
while the nonadaptive approach gives a bound in terms of the unknown
$p_i$s, but still gives useful asymptotics (Theorem~\ref{thm:whatV}).
The Bayesian solution 
sidesteps the various cases above, 
as it admits
a 
simple,
provably optimal empirical decision rule
(Section~\ref{sec:bayes}).
Unfortunately, 
we are unable to compute (or even nontrivially estimate)
the probability of error induced by this rule;
this is posed
as a challenging open problem.
\section{Background and related work}
Machine learning theory typically clusters {\em weighted majority}
\citep{DBLP:conf/focs/LittlestoneW89,DBLP:journals/iandc/LittlestoneW94}
within the framework of online algorithms;
see \citet{MR2409394} for a modern treatment.
Since the online setting is considerably more adversarial than ours,
we obtain very different 
weighted majority rules
and consistency guarantees.
The weights $w_i$ in (\ref{eq:N-P}) bear a striking similarity to the Adaboost update rule
\citep{261549,pre06049779}. However, the latter assumes weak learners with
access to labeled examples, while in our setting the experts are ``static''. Still, we do not
rule out a possible deeper connection between the Nitzan-Paroush decision rule and boosting.
In a recent line of work 
\citet{DBLP:conf/nips/LacasseLMGU06,DBLP:journals/jmlr/LavioletteM07,DBLP:conf/icml/RoyLM11}
have developed a PAC-Bayesian theory for the majority vote of simple classifiers. 
This approach facilitates data-dependent bounds and is even flexible enough to capture
some simple dependencies among the classifiers --- though, again, the latter are {\em learners}
as opposed to our {\em experts}. Even more recently, experts with adversarial noise have been considered
\citep{mansour2013}.
More directly related to the present work are  the papers of
\citet{berend-paroush-98}, which characterizes the consistency of the simple majority rule,
and
\citet{MR981253,DBLP:journals/scw/BerendS07}
which analyze various models of dependence among the experts.
\section{Known competences}
\label{sec:known}
In this section we 
assume that the expert competences $p_i$ are known
and
analyze the consistency of the Nitzan-Paroush optimal decision rule (\ref{eq:N-P}).
Our main result here is that the probability of error 
$\P(f\opt(\vec X)\neq Y)$
is small if and only if 
the committee potential $\Phi$ is large.
\hide{}
\begin{theorem}
\label{thm:N-P-known}
Suppose that the experts $\vec X=(X_1,\ldots,X_n)$ satisfy Assumptions (i)-(ii)
and $f:\set{\pm1}^n\to\set{\pm1}$ is the Nitzan-Paroush optimal decision rule.
Then
\begin{itemize}
\item[(i)]
${\ds
\P(f\opt(\vec X)\neq Y) \le \exp\paren{-\hf
\Phi
}
}. 
$
\item[(ii)] %
${\ds
\P(f\opt(\vec X)\neq Y) \ge
\frac{3}{4[1+\exp(
2\Phi+4\sqrt{\Phi}
)]}
}.
$
\end{itemize}
\end{theorem}
\hide{}
\paragraph{Open problem.} Exhibit (if possible) a function 
$g:\R\to\R$
such that
\shortlong{
$\P(f\opt(\vec X)\neq Y) \asymp g(\Phi)$.
}{
\beq
\P(f\opt(\vec X)\neq Y) &\asymp& g(\Phi).
\eeq
}
The remainder of this section
is devoted to proving Theorem~\ref{thm:N-P-known}.
\subsection{Proof of Theorem~\ref{thm:N-P-known}(i)}
\label{sec:pf(i)}
Define the $\set{0,1}$-indicator variables
\beqn
\label{eq:xidef}
\xi_i = \pred{X_i=Y},
\eeqn
corresponding to the event that the $i$\th expert is correct.
A mistake $f\opt(\vec X)\neq Y$ occurs precisely when\footnote{
Without loss of generality, ties are considered to be errors.}
the sum of the correct experts' weights fails to exceed
half the total mass:
\beqn
\label{eq:perrw2}
\P(f\opt(\vec X)\neq Y) = \P\paren{\sum_{i=1}^n w_i\xi_i \le \oo2\sum_{i=1}^n w_i}.
\eeqn
Since $\E\xi_i=p_i$, we may rewrite the 
probability in (\ref{eq:perrw2}) as
\shortlong{
\begin{align}
\label{eq:dev-rewrite}
& \P\paren{\sum_i w_i\xi_i \le \E\sqprn{\sum_i w_i\xi_i}-\sum_i (p_i-\hf)w_i}.
\end{align}
}{
\beqn
\label{eq:dev-rewrite}
\P\paren{\sum_i w_i\xi_i \le \E\sqprn{\sum_i w_i\xi_i}-\sum_i (p_i-\hf)w_i}.
\eeqn
}A standard tool for estimating such sum deviation probabilities is 
Hoeffding's inequality. 
Applied to (\ref{eq:dev-rewrite}), it yields the bound
\beqn
\label{eq:hoef-known}
\P(f\opt(\vec X)\neq Y) 
\le
\exp\paren{-\frac{2\sqprn{\sum_i (p_i-\hf)w_i}^2}{\sum_i w_i^2}},
\eeqn
which is far too crude for our purposes. 
Indeed, consider a finite committee of highly competent experts
with $p_i$'s arbitrarily close to $1$ and $X_1$ the most competent of all.
Raising $X_1$'s
competence
sufficiently far above his peers will cause both 
the numerator and the denominator in the exponent to be dominated
by $w_1^2$, making the right-hand-side of 
(\ref{eq:hoef-known}) bounded away from zero.
The inability of Hoeffding's inequality to guarantee consistency
even in such a felicitous setting is an instance of
its generally poor applicability to
highly heterogeneous sums,
a phenomenon explored in some depth in 
\citet{DBLP:journals/jmlr/McAllesterO03}.
Bernstein's and Bennett's inequalities suffer from a similar weakness
(see ibid.).
Fortunately, an inequality of
\citet{DBLP:conf/uai/KearnsS98} is sufficiently sharp
to yield the desired estimate: 
For all $p\in[0,1]$ and all $t\in\R$,
\beqn
\label{eq:k-s}
(1-p)e^{-tp}+pe^{t(1-p)}\le \exp\paren{\frac{1-2p}{4\log((1-p)/p)}t^2}
.
\eeqn
\begin{remark}
The Kearns-Saul inequality
(\ref{eq:k-s}) may be seen as a distribution-dependent
refinement of Hoeffding's (which
bounds the left-hand-side of (\ref{eq:k-s}) by $e^{t^2/8}$),
and
is not nearly as straightforward to prove. 
An elementary rigorous proof
is given in \citet{ECP2359}.
Following up,
\citet{DBLP:journals/corr/abs-1212-4663}
gave a ``soft'' proof
based on transportation and information-theoretic techniques.
\end{remark}
Put $\theta_i=\xi_i-p_i$, substitute into (\ref{eq:perrw2}),
and apply Markov's inequality:
\beqn
\label{eq:wtheta}
\P(f\opt(\vec X) \neq  Y) &=& 
\P\paren{ - \sum_{i} w_i\theta_i  \ge
\Phi
}\\
&\le&
\nonumber
e^{ -t \Phi }
 \E{\exp\paren{-t\sum_i w_i\theta_i}}.
\eeqn
Now
\beqn
\E e^{-t w_i\theta_i} &=& p_ie^{-(1-p_i)w_it}+(1-p_i)e^{p_iw_it} \nonumber\\
&\le& 
\exp\paren{\frac{-1+2p_i}{4\log(p_i/(1-p_i))}w_i^2t^2} 
\label{eq:use-KS}
\\\nonumber
&=& \exp\sqprn{\hf(p_i-\hf)w_it^2},
\eeqn
where the inequality follows from (\ref{eq:k-s}).
By independence,
\shortlong{
\begin{align*}
&\E \exp\paren{-t\sum_i w_i\theta_i} = \prod_i \E e^{-t w_i\theta_i} \\
&\le \exp\paren{\hf\sum_i(p_i-\hf)w_it^2}= \exp\paren{\hf \Phi t^2}
\end{align*}
}{
\beq
\E \exp\paren{-t\sum_i w_i\theta_i} &=& \prod_i \E e^{-t w_i\theta_i} \\
&\le& \exp\paren{\hf\sum_i(p_i-\hf)w_it^2}\\
&=& \exp\paren{\hf \Phi t^2}
\eeq
}and hence
\shortlong{
$\P(f\opt(\vec X)\neq Y) \le\exp\paren{\hf \Phi t^2-\Phi t}.$
}{
\beq
\P(f\opt(\vec X)\neq Y) &\le&\exp\paren{\hf \Phi t^2-\Phi t}.
\eeq
}Choosing $t=1$ yields the bound in Theorem~\ref{thm:N-P-known}(i).
\subsection{Proof of Theorem~\ref{thm:N-P-known}(ii)}
Define the $\set{\pm1}$-indicator variables
\beqn
\label{eq:etadef}
\eta_i = 2\pred{X_i=Y}-1,
\eeqn
corresponding to the event that the $i$\th expert is correct
and put $q_i=1-p_i$.
The shorthand $
\weta
=\sum_{i=1}^n w_i\eta_i$
will be convenient.
We will need some simple lemmata:
\begin{lemma}
\label{lem:antipodes}
\beq
\P(f\opt(\vec X)=Y) = \sum_{\vec\eta\in\set{\pm1}^n} \max\set{P(\vec\eta),P(-\vec\eta)}
\eeq
and
\beq
\P(f\opt(\vec X)\neq Y) = \sum_{\vec\eta\in\set{\pm1}^n} \min\set{P(\vec\eta),P(-\vec\eta)},
\eeq
where
\shortlong{
$P(\vec \eta) = 
\prod_{i:\eta_i=1}{p_i}
\prod_{i:\eta_i=-1}{q_i}
.
$
}{
\beq
P(\vec \eta) = 
\prod_{i:\eta_i=1}{p_i}
\prod_{i:\eta_i=-1}{q_i}
.
\eeq}
\end{lemma}
\bepf
The 
identities (\ref{eq:xidef}), (\ref{eq:perrw2}) and (\ref{eq:etadef})
imply 
that
a mistake occurs precisely when
$$ \sum_{i=1}^n w_i\frac{\eta_i+1}2\le\oo2\sum_{i=1}^n w_i,$$
which is equivalent to
\beqn
\label{eq:weta-err}
\weta
\le0.
\eeqn
Exponentiating both sides,
\shortlong{
\begin{align}
\exp\paren{
\weta
} 
&= \prod_{i=1}^n e^{w_i\eta_i}
\CC =\CC \CC 
\prod_{i:\eta_i=1}\frac{p_i}{q_i}
\cdot\CC \CC 
\prod_{i:\eta_i=-1}\frac{q_i}{p_i} 
=
\label{eq:h/h}
\frac{P(\vec\eta)}{P(-\vec\eta)} \le1.
\end{align}
}{\beqn
\exp\paren{
\weta
} 
&=& \prod_{i=1}^n e^{w_i\eta_i}
\nonumber\\
&=&
\prod_{i:\eta_i=1}\frac{p_i}{q_i}
\cdot
\prod_{i:\eta_i=-1}\frac{q_i}{p_i} \nonumber\\
&=&
\label{eq:h/h}
\frac{P(\vec\eta)}{P(-\vec\eta)} \le1.
\eeqn}We conclude from (\ref{eq:h/h}) that among two
``antipodal'' atoms $\pm\vec\eta\in\set{\pm1}^n$,
the one with the greater mass contributes to the probability being correct
and the one with the smaller mass contributes to the probability of error,
which proves the claim.
\enpf
\begin{remark}
The proof of Lemma~\ref{lem:antipodes}
also establishes
the optimality of the Nitzan-Paroush decision rule.
\end{remark}
\begin{lemma}
\label{lem:SS'}
Suppose that $\vec s,\vec s'\in(0,\infty)^m$
satisfy
\shortlong{
$\sum_{i=1}^m (s_i+s_i')\ge a$
}{
\beq
\sum_{i=1}^m (s_i+s_i')\ge a
\eeq}and
\shortlong{
$R\inv \le {s_i}/{s_i'} \le R$,
$i\in[m]$,
}{
\beq
\oo R \le \frac{s_i}{s_i'} \le R,
\qquad i\in[m]
\eeq}for some $R<\infty$.
Then
\shortlong{
$\sum_{i=1}^m \min\set{s_i,s_i'}\ge {a}/(1+R)$.
}{
\beq
\sum_{i=1}^m \min\set{s_i,s_i'}\ge\frac{a}{1+R}.
\eeq}
\end{lemma}
\bepf
Immediate from
\shortlong{
$s_i+s_i' \le \min\set{s_i,s_i'}(1+R)$.
}{
\beq
s_i+s_i' \le \min\set{s_i,s_i'}(1+R).
\eeq}
\enpf
\begin{lemma}
\label{lem:pqw}
Define the function $F:(0,1)\to\R$ by
\beq
F(x) = 
\frac{ x(1-x)\log(x/(1-x)) }{2x-1}.
\eeq
Then $\sup_{0<x<1}F(x)=\hf$.
\end{lemma}
\bepf
Deferred to the Appendix.
\enpf
Continuing with the main proof, observe that
\beqn
\label{eq:Eweta}
\E\sqprn{\weta} = \sum_{i=1}^n(p_i-q_i)w_i
= 2\Phi
\eeqn
and
\beq
\V\sqprn{\weta} = 4\sum_{i=1}^n p_iq_iw_i^2.
\eeq
By Lemma~\ref{lem:pqw},
\beq
p_iq_i w_i^2 \le \hf(p_i-q_i)w_i,
\eeq
and hence
\beqn
\label{eq:Vweta}
\V\sqprn{\weta} &\le& 4\Phi.
\eeqn
Define the segment $I\subset\R$ by
\beqn
\label{eq:Idef}
I=\left[
2\Phi-4\sqrt{\Phi},
2\Phi+4\sqrt{\Phi}
\right].
\eeqn
{
}Chebyshev's inequality together with (\ref{eq:Eweta}) and (\ref{eq:Vweta}) implies that
\beqn
\label{eq:3/4}
\P\paren{\weta\in I} &\ge& \frac34.
\eeqn
\hide{}
\hide{}
Consider an atom $\vec\eta\in\set{\pm1}^n$ for which
$\weta\in I$.
The proof of Lemma~\ref{lem:antipodes} shows that
\beqn
\frac{P(\vec\eta)}{P(-\vec\eta)} 
=
\exp\paren{
\weta
}
\le
\exp(
2\Phi+4\sqrt{\Phi}
)
,
\label{eq:hhW}
\eeqn
where the inequality follows from (\ref{eq:Idef}).
Lemma~\ref{lem:antipodes} further implies that
\beq
\P(f\opt(\vec X)\neq Y) &\ge& 
\sum_{
{\vec\eta\in\set{\pm1}^n:
\weta\in I
}} \min\set{P(\vec\eta),P(-\vec\eta)} \\
&\ge& \frac{3/4}{1+\exp(
2\Phi+4\sqrt{\Phi}
)},
\eeq
where the second inequality follows from Lemma~\ref{lem:SS'},
(\ref{eq:3/4}) and (\ref{eq:hhW}).
This completes the proof.
\section{Unknown competences: frequentist 
\shortlong{}{approach}
}
\label{sec:freq}
Our goal in this section is to obtain, insofar as possible,
analogues of Theorem~\ref{thm:N-P-known} for unknown expert competences.
When the 
${p_i}$s
are unknown, they must be estimated empirically before any
useful
weighted majority vote can be applied. 
There are various ways to model partial knowledge of expert competences \citep{
DBLP:journals/aamas/BaharadGKN11,
baharad12}.
Perhaps the simplest 
scenario
for
estimating the 
$p_i$s
is to assume that the $i$\th expert has been
queried 
independently
$m_i$ times, out of which he gave the correct prediction $k_i$ times.
Taking the $\set{m_i}$ to be fixed, define
the {\em committee profile} 
by
$\vec k=(k_1,\ldots,k_n)$;
this
is the aggregate of 
the agent's 
empirical knowledge of the experts' performance.
An {\em empirical decision rule}
$\hat f:
(\vec x,
\vec k
)\mapsto 
\set{\pm1}
$
makes a final decision based on the expert inputs $\vec x$
together with the committee profile.
Analogously to (\ref{eq:prerr}),
the probability of a mistake
is
\beqn
\label{eq:prerr-unk}
\P(\hat f(\vec X,\vec K)\neq Y).
\eeqn
Note that now the committee profile is an additional source of randomness.
Here we run into our first difficulty: unlike 
the probability
in (\ref{eq:prerr}),
which is minimized by the Nitzan-Paroush rule, 
the agent cannot formulate
an optimal decision rule $\hat f$ in advance
without knowing the ${p_i}$s.
This is because
no decision rule is optimal uniformly over the range
of possible $p_i$s.
\hide{}
Our approach will be to consider weighted majority decision rules of the form
\beqn
\label{eq:majvote-freq}
\hat f(\vec x,\vec k) = \sgn\paren{\sum_{i=1}^n \hat w(k_i) x_i}
\eeqn
and to analyze their consistency properties under two different regimes:
low-confidence and high-confidence.
These refer to the confidence intervals of the frequentist estimate of $p_i$, given by
\beqn
\label{eq:pihat}
\hat p_i = \frac{k_i}{m_i}.
\eeqn
\subsection{Low-confidence regime}
In the low-confidence regime, the sample sizes $m_i$ may be as small as $1$,
and we define\footnote{
For 
$m_i\min\set{p_i,q_i}\ll1$,
the estimated competences $\hat p_i$ 
may well
take values in $\set{0,1}$,
in which case
$\log(\hat p_i/\hat q_i)={\pm\infty}$.
The rule in
(\ref{eq:what-freq})
is essentially a 
first-order Taylor approximation to $w(\cdot)$ about 
$p=\hf$.
}
\beqn
\label{eq:what-freq}
\hat w(k_i) = \hat w_i\lc := \hat p_i-\hf, \qquad i\in[n],
\eeqn
which induces the empirical decision rule $\hat f\lc$.
It remains to analyze $\hat f\lc$'s probability of error. Recall the definition of $\xi_i$ from
(\ref{eq:xidef}) and observe that
\beqn
\label{eq:Ewxi}
\E
\sqprn{
\hat w_i\lc\xi_i
}
=
\E[ (\hat p_i-\hf)\xi_i]
=
( p_i-\hf)p_i,
\eeqn
since $\hat p_i$ and $\xi_i$ are independent.
As in (\ref{eq:perrw2}), the probability of error
(\ref{eq:prerr-unk}) is 
\beqn
\label{eq:PZ}
\P\paren{\sum_{i=1}^n \hat w_i\lc\xi_i \le \oo2\sum_{i=1}^n \hat w_i\lc}
=
\P\paren{\sum_{i=1}^n Z_i\le0},
\eeqn
where $Z_i=\hat w_i\lc(\xi_i-\hf)$. Now the $\set{Z_i}$ are independent
random variables,
$\E Z_i=(p_i-\hf)^2$
(by (\ref{eq:Ewxi})),
and each $Z_i$ takes values in an interval of length $\hf$.
Hence, the standard Hoeffding bound applies:
\beqn
\label{eq:low-hoeff}
\P(\hat f\lc(\vec X,\vec K)\neq Y)
\le
\exp\sqprn{-\frac8n\paren{\sum_{i=1}^n(p_i-\hf)^2}^2}.
\eeqn
We summarize these calculations in 
\begin{theorem}
\label{thm:unknown-freq-low}
A sufficient condition for 
$\P(\hat f\lc(\vec X,\vec K)\neq Y)\to0$
is
\shortlong{
$\oo{\sqrt n}\sum_{i=1}^n(p_i-\hf)^2
\to\infty
$.
}{
\beq
\oo{\sqrt n}\sum_{i=1}^n(p_i-\hf)^2
\to\infty
.
\eeq}
\end{theorem}
Several remarks are in order.
First, notice that the error bound in (\ref{eq:low-hoeff})
is stated in terms of the unknown $\set{p_i}$,
providing the agent with 
large-committee
asymptotics but 
giving no
finitary
information;
this limitation is inherent in the low-confidence regime.
Secondly, 
the condition in
Theorem~\ref{thm:unknown-freq-low} is considerably
more restrictive
than 
the consistency condition
$\Phi\to\infty$
implicit in Theorem~\ref{thm:N-P-known}.
Indeed, the empirical decision rule
$\hat f\lc$
is incapable of exploiting a single highly competent expert
in the way that 
$f\opt$ from
(\ref{eq:N-P}) does. 
Our analysis could be sharpened somewhat for moderate sample sizes $\set{m_i}$
by using Bernstein's inequality to take advantage of the low variance of the $\hat p_i$s.
For sufficiently large sample sizes, however, the high-confidence regime (discussed below) begins to take over.
Finally, 
there is one sense in which
this case is ``easier'' to analyze than that of known $\set{p_i}$:
since the summands in (\ref{eq:PZ}) are bounded, Hoeffding's inequality
gives nontrivial results and there is no need for more advanced tools
such as the Kearns-Saul inequality (\ref{eq:k-s}) (which is actually inapplicable in this case).
\subsection{High-confidence regime}
\newcommand{\tips}{\tilde\eps}
\newcommand{\sqg}{\sqrt\gamma}
\newcommand{\sqgi}{\sqrt{\gamma_i}}
In the high-confidence regime,
each estimated competence $\hat p_i$
is close to the true value $p_i$ with high probability.
To formalize this, 
fix some $0<\delta<1$, $0<\eps\le5$, 
and 
put
\beq
q_i = 1 - p_i,
~
\hat q_i = 1 - \hat p_i
.
\eeq
We will set
the empirical weights according to
the ``plug-in'' Nitzan-Paroush rule
\beqn
\label{eq:w-plug}
\hat w_i\hc := \log\frac{\hat p_i}{\hat q_i},
\qquad i\in[n]
,
\eeqn
which 
induces the empirical decision rule $\hat f\hc$
and raises immediate concerns about
$\hat w_i\hc=\pm\infty$.
We give two kinds of bounds on
$\P(\hat f\hc\neq Y)$:
nonadaptive
and adaptive.
In the nonadaptive analysis,
we show that for $m_i\minpq_i\gg1$,
with high probability $\abs{w_i-\hat w_i\hc}\ll1$,
and thus a ``perturbed'' version of 
Theorem~\ref{thm:N-P-known}(i) holds
(and in particular, $w_i\hc$ will be finite with high probability).
In the adaptive analysis, we allow
$\hat w_i\hc$ to take on infinite values\footnote{
When the decision rule is faced with evaluating sums involving $\infty-\infty$,
we automatically count this as an error.
} and show
(perhaps surprisingly) that this decision rule
also asymptotically achieves the 
rate of
Theorem~\ref{thm:N-P-known}(i).
\paragraph{Nonadaptive analysis.}
Define $\tips\in(0,1)$ by $\eps=2\tips+4\tips^2$
or, explicitly,
\beqn
\label{eq:tips}
\hide{}
\tips = \frac{\sqrt{4\eps+1}-1}4.
\eeqn
\begin{lemma}
\label{lem:emtau}
If
\beqn
\label{eq:emtau}
\tips^2 
m_ip_i &\ge& 3\log({2n}/{\delta})
,\qquad i\in[n],
\eeqn
then
\shortlong{
$
\P\paren{
\exists i\in[n]:
\frac{\hat p_i}{p_i}\notin
(1-\tips,1+\tips)
}
\le\delta
.
$
}{
\beq
\P\paren{
\exists i\in[n]:
\frac{\hat p_i}{p_i}\notin
(1-\tips,1+\tips)
}
\le\delta
.
\eeq}
\end{lemma}
\bepf
The multiplicative Chernoff bound yields
\beq
\P\paren{\hat p_i < (1-\tips)p_i} \le e^{-\tips^2m_ip_i/2}
\eeq
and
\beq
\P\paren{\hat p_i > (1+\tips)p_i} \le e^{-\tips^2m_ip_i/3}.
\eeq
Hence,
\shortlong{
$
\P\paren{
\frac{\hat p_i}{p_i}
\notin
(1-\tips,1+\tips)
}
\le
2e^{-\tips^2m_ip_i/3}.
$
}{
\beq
\P\paren{
\frac{\hat p_i}{p_i}
\notin
(1-\tips,1+\tips)
}
&\le&
2e^{-\tips^2m_ip_i/3}.
\eeq}
The claim follows from (\ref{eq:emtau}) and the union bound.
\enpf
\hide{}
\begin{lemma}
\label{lem:wwee}
Let $w_i$ be the optimal Nitzan-Paroush weight (\ref{eq:wdef}).
If
\shortlong{
$1-\tips\le {\hat p_i}/{p_i}, {\hat q_i}/{q_i}\le1+\tips$,
then $\abs{ w_i-\hat w_i\hc}\le \eps$.
}{
\beq
1-\tips\le
\frac{\hat p_i}{p_i}
,
\frac{\hat q_i}{q_i}
\le1+\tips
\eeq
then
\beq
\abs{ w_i-\hat w_i\hc}\le \eps.
\eeq}
\end{lemma}
\bepf
We have
\shortlong{
\begin{align*}
\abs{
w_i-\hat w_i\hc
}
=&
\abs{
\log\frac{p_i}{q_i}
-
\log\frac{\hat p_i}{\hat q_i}
} 
=
\abs{
\log\frac{p_i}{\hat p_i}
-
\log\frac{\hat q_i}{ q_i}
}\\
=&
\abs{
\log\frac{p_i}{\hat p_i}
}
+
\abs{
\log\frac{\hat q_i}{ q_i}
}.
\end{align*}
}{
\beq
\abs{
w_i-\hat w_i\hc
}
&=&
\abs{
\log\frac{p_i}{q_i}
-
\log\frac{\hat p_i}{\hat q_i}
} \\
&=&
\abs{
\log\frac{p_i}{\hat p_i}
-
\log\frac{\hat q_i}{ q_i}
}\\
&=&
\abs{
\log\frac{p_i}{\hat p_i}
}
+
\abs{
\log\frac{\hat q_i}{ q_i}
}.
\eeq}Now\footnote{
The first containment requires
$\log(1-x)\ge -x-2x^2$,
which holds 
(not exclusively)
on 
$(0,0.9)$.
The restriction $\eps\le5$ ensures that $\tips$ is in this range.
}\shortlong{
\beq
[ \log(1-\tips) ,\log(1+\tips) ]
\subseteq 
[ -\tips-2\tips^2 ,\tips ] 
\subseteq
[-\hf\eps,\hf\eps]
,
\eeq
whence
$
\abs{
\log\frac{p_i}{\hat p_i}
}
+
\abs{
\log\frac{\hat q_i}{ q_i}
}
\le\eps.
$
}{
\beq
[ \log(1-\tips) ,\log(1+\tips) ]
&\subseteq & 
[ -\tips-2\tips^2 ,\tips ] 
\\
&\subseteq & 
[-\hf\eps,\hf\eps]
,
\eeq
whence
\beq
\abs{
\log\frac{p_i}{\hat p_i}
}
+
\abs{
\log\frac{\hat q_i}{ q_i}
}
\le\eps.
\eeq
}
\enpf
\begin{corollary}
\label{cor:|ww|}
\shortlong{
If
$
\tips^2 m_i\minpq_i \ge 3\log({4n}/{\delta})
$,
$i\in[n]$,
then
$
\P\paren{
\max_{i\in[n]}
\abs{
w_i-\hat w_i\hc
} >\eps} \le\delta.
$
}{
If
\beq
\tips^2 m_i\minpq_i &\ge& 3\log({4n}/{\delta}),
\qquad
i\in[n],
\eeq
then
\beq
\P\paren{
\max_{i\in[n]}
\abs{
w_i-\hat w_i\hc
} >\eps} &\le&\delta.
\eeq}
\end{corollary}
\bepf
An immediate consequence of applying Lemma~\ref{lem:emtau}
to $p_i$ and $q_i$ with the union bound.
\enpf
\newcommand{\bw}{\vec w}
\newcommand{\bwhat}{\hat{\vec w}\hc}
To state the next result,
let us arrange
the plug-in weights (\ref{eq:w-plug})
as a vector $\bwhat\in\R^n$,
as was done with
$\bw$ and $\veta$ from 
Section~\ref{sec:pf(i)}.
The corresponding weighted majority rule $\hat f\hc$
yields
an error
precisely when
\shortlong{
$\bwhat \cdot \veta \le 0$
}{
\beq
\bwhat \cdot \veta \le 0
\eeq}(cf. (\ref{eq:weta-err})).
Our nonadaptive approach culminates in the following result.
\begin{theorem}
\label{thm:whatV}
Let 
$0<\delta<1$
and
$0<\eps<\min\set{5,2\Phi/n}$.
If
\beqn
\label{eq:thm:whatV-Assuming}
\CC\CC\CC
m_i\minpq_i \ge 
3\CC
\paren{\CC\frac{\sqrt{4\eps+1}-1}4\CC}^{\CC\CC-2}
\CC\CC\CC\CC\CC\CC\log\frac{4n}{\delta},
\CC\CC\CC\CC\CC\CC\CC
\qquad i\in[n],
\eeqn
then
\beqn
\label{eq:thm:whatV-we have}
\P\paren{
\hat f\hc(\vec X,\vec K)\neq Y
}
\;\le\;
\delta+
\exp\sqprn{
-\frac{ (2\Phi-\eps n)^2}{8\Phi}}.
\eeqn
\end{theorem}
\begin{remark}
For fixed $\set{p_i}$ and $\min_{i\in[n]}m_i\to\infty$,
we may take $\delta$ and $\eps$ arbitrarily small
--- and in this limiting case, the bound of
Theorem~\ref{thm:N-P-known}(i) is recovered.
\end{remark}
\bepf[Proof of Theorem~\ref{thm:whatV}]
Since
\beq
\abs{ 
\bw\cdot\veta
- 
\bwhat \cdot \veta 
} &=&
\abs{ (\bw-\bwhat)\cdot\veta} \\
&\le& \sum_{i=1}^n\abs{w_i-w_i\hc} 
=\nrm{\bw-\bwhat}_1,
\eeq
we have
\beq
\P\paren{\bwhat \CC\cdot\CC \veta \le 0}
\le 
\P(\nrm{\bw-\bwhat}_1>\eps n)
+
\P(\bw\CC\cdot\CC\veta\le\eps n).
\eeq
Corollary~\ref{cor:|ww|} upper-bounds
the first term on the right-hand side 
by $\delta$. The second term
is estimated by
replacing 
$\Phi$
by
$
\Phi
-\eps n$
in (\ref{eq:wtheta})
and repeating the argument
following that formula.
\enpf
\paragraph{Adaptive analysis.}
Theorem~\ref{thm:whatV} has the drawback of being {\em nonadaptive},
in that
its
assumptions
(\ref{eq:thm:whatV-Assuming})
and 
conclusions
(\ref{eq:thm:whatV-we have})
depend on the unknown $\set{p_i}$
and hence cannot be evaluated by the agent (the 
bound in (\ref{eq:low-hoeff})
is also nonadaptive).
In the {\em adaptive} approach, all results are stated
in terms of empirically observed quantities:
\begin{theorem}
\label{thm:adapt}
Put\footnote{
\label{ftn:inex}
Actually, as the proof will show,
we may take $\delta$ to be a smaller value, but with a more
complex dependence on $\set{m_i}$, which simplifies to
$2[1-(1-(2\sqrt m)\inv)^n]$ for $m_i\equiv m$.
}
\beq
\delta = \sum_{i=1}^n\oo{\sqrt{m_i}}
\eeq
and
let $R$
be the event
\beqn
\label{eq:R}
\exp\paren{-\oo2\sum_{i=1}^n(\hat p_i-\hf)\hat w_i\hc}
&\le& \frac{\delta}{2}.
\eeqn
Then
\beq
\P
\paren{R\LAND 
\set{\hat f\hc(\vec X,\vec K)\neq Y}}
&\le &
\delta.
\eeq
\end{theorem}
\begin{remark}
\label{rem:adapt}
Our 
interpretation for Theorem~\ref{thm:adapt} is as follows.
The agent 
observes
the committee profile $\vec K$, 
which determines
the 
$\set{\hat p_i,\hat w_i\hc}$,
and
then checks whether
the event $R$ has occurred.
If not, the adaptive agent refrains from
making a decision 
(and may choose to fall back
on the low-confidence approach described previously). 
If $R$
{\em does} hold, however, the agent predicts $Y$ according to $\hat f\hc$.
As explained above, 
there 
does not exist a
nontrivial a priori upper bound on
$\P(\hat f\hc(\vec X,\vec K)\neq Y)$
absent any knowledge of the ${p_i}$s.
Instead, Theorem~\ref{thm:adapt} bounds the probability of the agent being
``fooled'' by 
an unrepresentative
committee profile.\footnote{These adaptive bounds are similar in spirit to
{\em empirical Bernstein} methods, 
\citep{DBLP:conf/alt/AudibertMS07,DBLP:conf/icml/MnihSA08,DBLP:conf/colt/MaurerP09},
where the player's confidence depends on the empirical variance.
}
Observe that for 
$\min_{i\in[n]}m_i\gg1$,
we have $\hat p_i\approx p_i$
and $\hat w_i\hc\approx w_i$,
and thus Theorem~\ref{thm:N-P-known}(i)
is recovered as a limiting case.
Note that we have done nothing to prevent $\hat w_i\hc=\pm\infty$,
and this may indeed happen. Intuitively, 
there are two reasons for infinite
$\hat w_i\hc$:
(a) noisy $\hat p_i$ due to $m_i$ being too small, or
(b) the $i$\th expert is actually highly (in)competent, which causes
$\hat p_i\in\set{0,1}$ even for large $m_i$.
The $1/{\sqrt{m_i}}$ term in the bound ensures against case (a),
while 
in case (b), choosing 
infinite
$\hat w_i\hc$ causes no harm
(as we show in the proof).
\end{remark}
\bepf[Proof\shortlong{ of Theorem \ref{thm:adapt}}{}]
We will write the probability and expectation operators 
with subscripts (such as $\vec K$) to indicate the random variable(s) being summed over.
Thus,
\shortlong{
\begin{align}
\nonumber& 
\P_{\vec K,\vec X,Y} \paren{R \LAND 
\set{\hat f\hc(\vec X,\vec K)\neq Y} }
\\\nonumber=&
\P_{\vec K,\veta}\paren{R \LAND \set{\bwhat \cdot \veta \le 0}}
\\=&
\E_{\vec K}\sqprn{ \chr_{R} \cdot \P_{\veta}\paren{\bwhat \cdot \veta \le 0 \gn \vec K}}.
\label{eq:EP-decomp}
\end{align}
}{
\beqn
\nonumber
\P_{\vec K,\vec X,Y} \paren{R \LAND 
\set{\hat f\hc(\vec X,\vec K)\neq Y} }
&=&
\P_{\vec K,\veta}\paren{R \LAND 
\set{\bwhat \cdot \veta \le 0}
}\\
\nonumber
&=& \E_{\vec K}\sqprn{ \chr_{R} \cdot \P_{\veta}\paren{\bwhat \cdot \veta \le 0 \gn \vec K}}.
\\
\label{eq:EP-decomp}
\eeqn}Recall that
the random variable 
$\veta\in\set{\pm1}^n$,
with
probability mass function
\shortlong{
$
P(\veta) = 
\prod_{i:\eta_i=1}{p_i}
\prod_{i:\eta_i=-1}{q_i},
$
}{
\beq
P(\veta) = 
\prod_{i:\eta_i=1}{p_i}
\prod_{i:\eta_i=-1}{q_i},
\eeq}is independent of
$\vec K$,
and hence
\beqn
\label{eq:Pkh}
\P_{\veta}\paren{\bwhat \cdot \veta \le 0 \gn \vec K}
=
\P_{\veta} \paren{\bwhat \cdot \veta \le 0}.
\eeqn
Define the random variable
$\hat\veta\in\set{\pm1}^n$ 
(conditioned on $\vec K$)
by the probability mass function
\shortlong{
$
P(\hat\veta) = 
\prod_{i:\eta_i=1}{\hat p_i}
\prod_{i:\eta_i=-1}{\hat q_i},
$
}{
\beq
P(\hat\veta) = 
\prod_{i:\eta_i=1}{\hat p_i}
\prod_{i:\eta_i=-1}{\hat q_i},
\eeq}and the set $A\subseteq\set{\pm1}^n$ by
$
A=\set{\vec x : \bwhat \cdot \vec x\le0}.
$
Now
\shortlong{
\begin{align*}
&\abs{ 
\P_{\veta} \paren{\bwhat \cdot \veta \le 0}-
\P_{\hat\veta} \paren{\bwhat \cdot \hat\veta \le 0}}\\
&=
\abs{\P_{\veta} \paren{A}-
\P_{\hat\veta} \paren{A}
}
\le
\max_{A\subseteq\set{\pm1}^n}\abs{
\P_{\veta} \paren{A}-
\P_{\hat\veta} \paren{A}
}\\
&=\TV{\P_{\veta} -\P_{\hat\veta}}
\le \sum_{i=1}^n \abs{p_i-\hat p_i}=:M
,
\end{align*}
}{
\beq
\nonumber
\abs{ 
\P_{\veta} \paren{\bwhat \cdot \veta \le 0}-
\P_{\hat\veta} \paren{\bwhat \cdot \hat\veta \le 0}}
&=&
\nonumber
\abs{\P_{\veta} \paren{A}-
\P_{\hat\veta} \paren{A}
}\\
&\le&
\nonumber
\max_{A\subseteq\set{\pm1}^n}\abs{
\P_{\veta} \paren{A}-
\P_{\hat\veta} \paren{A}
}\\
\nonumber
&=&\TV{\P_{\veta} -\P_{\hat\veta}} \\
\nonumber
\label{eq:alpha-sum}
&\le& \sum_{i=1}^n \abs{p_i-\hat p_i}
=:M
,
\eeq}where 
the 
inequality follows from a standard
tensorization
property
of the total variation norm $\TV{\cdot}$,
see e.g. \citep[Lemma 2.2]{kontorovich12}.\hide{}
By Theorem~\ref{thm:N-P-known}(i), we have
\beq
\P_{\hat\veta} \paren{\bwhat \cdot \hat\veta \le 0}
\le \exp\paren{-\hf\sum_{i=1}^n(\hat p_i-\hf)\hat w_i\hc},
\eeq
and hence
\beq
\P_{\veta} \paren{\bwhat \cdot \veta \le 0}
\le 
M+
\exp\paren{-\hf\sum_{i=1}^n(\hat p_i-\hf)\hat w_i\hc}.
\eeq
Invoking (\ref{eq:Pkh}), we substitute the right-hand side above
into (\ref{eq:EP-decomp}) to obtain
\shortlong{
\begin{align*}
&\P_{\vec K,\vec X,Y} \paren{R \LAND 
\set{\hat f\hc(\vec X,\vec K)\neq Y} }
\\&\le
\E_{\vec K}\sqprn{ \chr_{R} \cdot 
\paren{
M+
\exp\paren{-\hf\sum_{i=1}^n(\hat p_i-\hf)\hat w_i\hc}
}
}
\\&\le
\E_{\vec K}[M]
+
\E_{\vec K}\sqprn{
\chr_{R}
\exp\paren{-\hf\sum_{i=1}^n(\hat p_i-\hf)\hat w_i\hc}
}.
\end{align*}
}{
\beq
\P_{\vec K,\vec X,Y} \paren{R \LAND 
\set{\hat f\hc(\vec X,\vec K)\neq Y} }
&\le&
\E_{\vec K}\sqprn{ \chr_{R} \cdot 
\paren{
M+
\exp\paren{-\hf\sum_{i=1}^n(\hat p_i-\hf)\hat w_i\hc}
}
}\\
&\le&
\E_{\vec K}[M]
+
\E_{\vec K}\sqprn{
\chr_{R}
\exp\paren{-\hf\sum_{i=1}^n(\hat p_i-\hf)\hat w_i\hc}
}.
\eeq}By 
the
definition of $R$, the second term on the last right-hand side
is 
upper-bounded by $\delta/2$.
To estimate $M$, we invoke a simple mean absolute deviation bound
(cf. \citet{Berend2013}):
\beq
\E_{\vec K}
\abs{p_i-\hat p_i} \;\le\; \sqrt{\frac{p_i(1-p_i)}{m_i}} \;\le\; \oo{2\sqrt{m_i}},\\
\eeq
which finishes the proof.
\enpf
\begin{remark}
The improvement mentioned in Footnote~\ref{ftn:inex}
is achieved
via a refinement of the bound
$\TV{\P_{\veta} -\P_{\hat\veta}}\le \sum_{i=1}^n \abs{p_i-\hat p_i}$
to
$\TV{\P_{\veta} -\P_{\hat\veta}}\le \alpha\paren{\set{ \abs{p_i-\hat p_i} : i\in[n]}}$,
where
$\alpha(\cdot)$ is the 
function defined in
\citet[Lemma 4.2]{kontorovich12}.
\end{remark}
\paragraph{Open problem.} 
As argued in Remark~\ref{rem:adapt},
Theorem~\ref{thm:adapt} achieves
the optimal asymptotic rate in $\set{p_i}$.
Can the dependence
on $\set{m_i}$ be improved,
perhaps through a better choice of $\bwhat$?
\section{Unknown competences: Bayesian
\shortlong{}{approach}
}
\label{sec:bayes}
A shortcoming of Theorem~\ref{thm:adapt}
is that 
when condition $R$ fails, the agent
is left with no estimate of the error probability.
\hide{}An alternative (and in some sense cleaner) approach to handling
unknown expert competences $p_i$ is 
to assume a known prior distribution over the competence levels $p_i$.
The natural choice of prior for a Bernoulli parameter is the Beta distribution, namely
\shortlong{
$p_i \sim \Bdist(\alpha_i,\beta_i)$
}{
\beq
p_i \sim \Bdist(\alpha_i,\beta_i)
\eeq}with density
\beq
\frac{p_i^{\alpha_i-1}q_i^{\beta_i-1}}{B(\alpha_i,\beta_i)},
\qquad \alpha_i,\beta_i>0,
\eeq
where $q_i=1-p_i$ and $B(x,y)=\Gamma(x)\Gamma(y)/\Gamma(x+y)$.
Our full probabilistic model is as follows.
Each of the $n$ expert competences $p_i$ is drawn independently
from a Beta distribution with known parameters $\alpha_i,\beta_i$.
Then the $i$\th expert, $i\in[n]$,
is
queried 
independently
$m_i$ times, 
with $k_i$ correct predictions and $m_i-k_i$ incorrect ones. As before,
$\vec K=(k_1,\ldots,k_n)$ 
is
the
(random)
{committee profile}.
Absent direct knowledge of the $p_i$s, the agent relies on
an {empirical decision rule}
$\hat f:
(\vec x,
\vec k
)\mapsto 
\set{\pm1}
$
to produce a final decision based on the expert inputs $\vec x$
together with the committee profile $\vec k$.
A decision rule 
$\hat f\ba
$ 
is {\em Bayes-optimal} if it minimizes
\beqn
\label{eq:prerr-bayes}
\P(\hat f(\vec X,\vec K)\neq Y),
\eeqn
which is formally identical to 
(\ref{eq:prerr-unk})
but semantically there is a difference:
the probability in (\ref{eq:prerr-bayes}) is over the $p_i$ in addition
to $(\vec X,Y,\vec K)$.
Unlike the frequentist approach,
where no optimal empirical decision rule was possible,
the Bayesian approach readily admits one.
For a given $\vec x\in\set{\pm1}$, define $I_+(\vec x)$
to be the set of YES votes
\beq
I_+(\vec x) = \set{i\in [n]: x_i=+1}
\eeq
and $I_-(\vec x)=[n]\setminus I_+(\vec x)$ to be the set of NO votes.
Let us fix some $A\subseteq[n]$, $B=[n]\setminus A$ and compute
\begin{align*}
&\P(Y=+1, I_+(\vec X)=A, I_-(\vec X)=B)
\\
&=\prod_{i=1}^n \int_0^1
\frac{p_i^{\alpha_i-1}q_i^{\beta_i-1}}{B(\alpha_i,\beta_i)}
\binom{m_i}{k_i}p_i^{k_i}q_i^{m_i-k_i}
p_i^{\pred{i\in A}}q_i^{\pred{i\in B}}dp_i
\\
&=\prod_{i=1}^n \CC \CC 
\frac{ \binom{m_i}{k_i}p_i^{k_i}q_i^{m_i-k_i} } {B(\alpha_i,\beta_i)}
\CC \CC \CC \int_0^1 \CC \CC \CC \CC \CC  p_i^{\alpha_i\CC +k_i\CC -1+\pred{i\in A}}\CC q_i^{\beta_i\CC +m_i\CC -k_i\CC -1+\pred{i\in B}}\CC dp_i
\\
&=\prod_{i=1}^n 
\frac{
\binom{m_i}{k_i}B(\alpha_i+k_i+\pred{i\in A},\beta_i+m_i-k_i+\pred{i\in B})
}{B(\alpha_i,\beta_i)}
.
\end{align*}
Analogously,
\begin{align*}
&\P(Y=-1, I_+(\vec X)=A, I_-(\vec X)=B)
\\
&=\prod_{i=1}^n 
\frac{
\binom{m_i}{k_i}B(\alpha_i+k_i+\pred{i\in B},\beta_i+m_i-k_i+\pred{i\in A})
}{B(\alpha_i,\beta_i)}
.
\end{align*}
Let us use the shorthand $P(+1,A,B)$ and $P(-1,A,B)$ for the joint probabilities in the last two displays, 
along 
with their
corresponding conditionals $P(\pm 1\gn A,B)$. Obviously, 
\beq
P(1| A,B)\CC >\CC P(-1| A,B) \CC \CC \iff\CC \CC 
P(1, A,B)\CC >\CC P(-1, A,B),
\eeq
which occurs precisely if
\begin{align*}
 &\prod_{i=1}^n B(\alpha_i+k_i+\pred{i\in A},\beta_i+m_i-k_i+\pred{i\in B})\\
&>\prod_{i=1}^n B(\alpha_i+k_i+\pred{i\in B},\beta_i+m_i-k_i+\pred{i\in A}).
\end{align*}
This is equivalent to
\begin{align*}
 \prod_{i\in A} \CC (\alpha_i\CC +\CC k_i)\CC \prod_{i\in B}\CC (\beta_i\CC +\CC m_i\CC -\CC k_i)
\CC >\CC \prod_{i\in B} \CC (\alpha_i\CC +\CC k_i)\CC \prod_{i\in A}\CC (\beta_i\CC +\CC m_i\CC -\CC k_i)
,
\end{align*}
which further simplifies to
\beq
\prod_{i\in A}\frac{\alpha_i+k_i}{\beta_i+m_i-k_i}
>
\prod_{i\in B}\frac{\alpha_i+k_i}{\beta_i+m_i-k_i}.
\eeq
Hence, the choice
\beq
\hat w_i\ba =\log\frac{\alpha_i+k_i}{\beta_i+m_i-k_i}
\eeq
guarantees that the decision rule
\shortlong{
$
\hat f\ba(\vec x,\vec k) = \sgn\paren{\sum_{i=1}^n\hat w_i\ba x_i}
$
}{
\beq
\hat f\ba(\vec x,\vec k) = \sgn\paren{\sum_{i=1}^n\hat w_i\ba x_i}
\eeq}maximizes the probability of being correct for each input $\vec x\in\set{\pm1}^n$.
Notice that for $0<p_i<1$, we have
\shortlong{
$
\hat w_i\ba {\longto_{m_i\to\infty}}w_i,
$
}{
\beq
\hat w_i\ba 
{\longto_{m_i\to\infty}}
w_i,
\qquad i\in[n],
\eeq}
almost surely, both in the frequentist and the Bayesian interpretations.
Unfortunately, although
\beq
\P(\hat f\ba(\vec X,\vec K)\neq Y)
=
\P(
\hat{\vec w}\ba
\cdot\veta \le 0
)
\eeq
is a deterministic function of $\set{\alpha_i,\beta_i,m_i}$,
we are unable to compute it at this point, or even give a non-trivial bound.
One source of the problem is the 
coupling between
$\hat{\vec w}\ba$ and~$\veta$.
\paragraph{Open problem.} Give a non-trivial estimate for
$\P(\hat f\ba(\vec X,\vec K)\neq Y)$.
\section{Experiments}
\begin{figure}[ht]
\vskip 0.2in
\begin{center}
\centerline{\includegraphics[width=1.15\columnwidth]{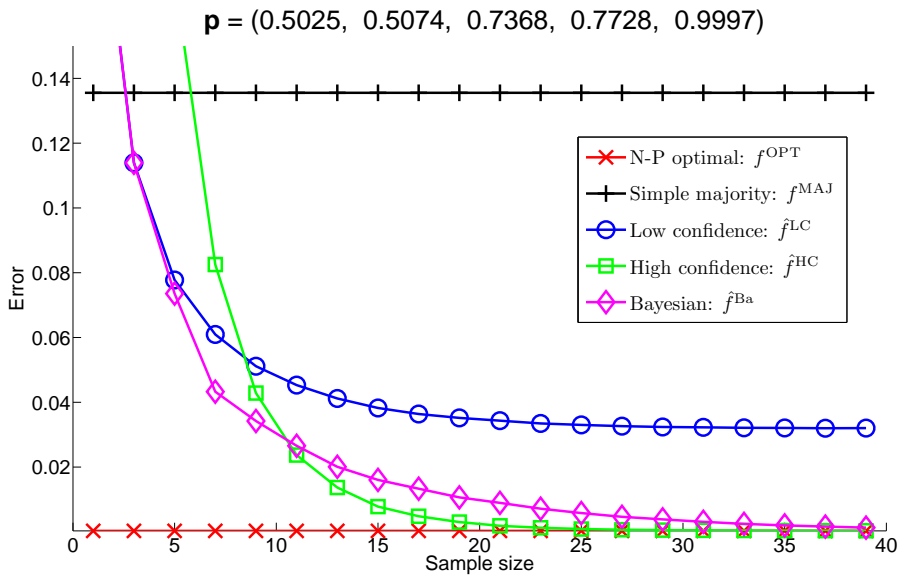}}
\caption{For very small sample sizes, $\hat f\lc$ outperforms
$\hat f\hc$ but is outperformed by $\hat f\ba$. Starting from sample size
$\approx13$, $\hat f\hc$ dominates the other empirical rules.
The empirical rules are (essentially) sandwiched between $f\opt$ and $f\maj$.
}
\label{fig:special5}
\end{center}
\vskip -0.2in
\end{figure}

\begin{figure}[H]
\vskip 0.2in
\begin{center}
\centerline{\includegraphics[width=1.15\columnwidth]{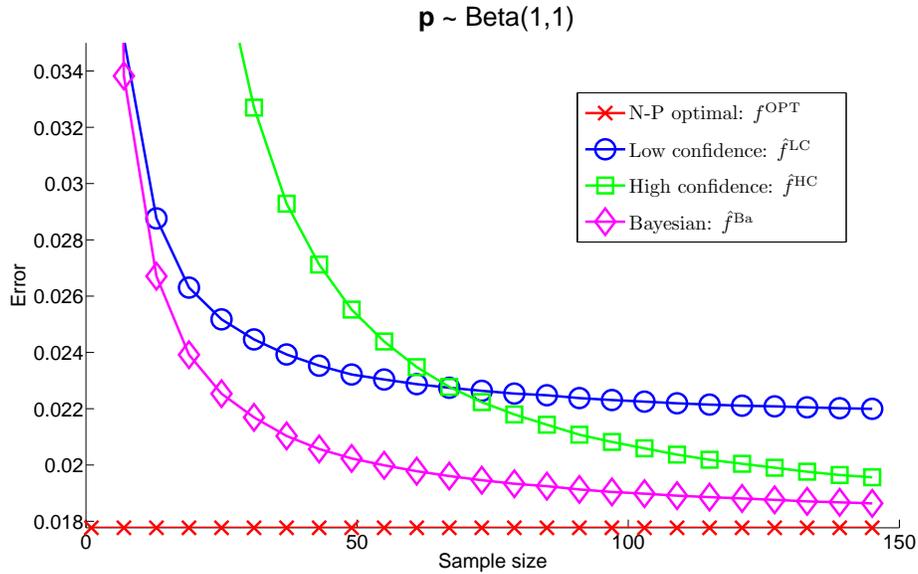}}
\caption{Unsurprisingly, $\hat f\ba$ uniformly outperforms the other
two empirical rules. 
We found it
somewhat surprising 
that
$\hat f\hc$ 
required so many samples (about $60$ on average)
to overtake $\hat f\lc$.
The simple majority rule $f\maj$ (off the chart) performed at 
an average accuracy of $50\%$, as expected.
}
\label{fig:bayes}
\end{center}
\vskip -0.2in
\end{figure}

It is most instructive to 
take the committee size $n$
to be small
when comparing the different voting rules.
Indeed, for a large committee of ``marginally competent'' experts
with $p_i=\hf+\gamma$ for some $\gamma>0$,
even
the simple majority rule
$f\maj(\vec x)=\sgn(\sum_{i=1}^n x_i)$
has a probability of error decaying as
$\exp(-
4n\gamma^2
)$,
as can be easily seen from Hoeffding's bounds.
The more sophisticated voting rules
discussed in this paper 
perform even better 
in this setting.
Hence, small committees provide the natural test-bed
for gauging a voting rule's ability to 
exploit highly competent experts.
In our experiments, we set $n=5$ and the sample sizes $m_i$ were identical
for all experts.
The results were averaged over $10^5$ trials.
Two of our experiments are described below.
\paragraph{Low vs. high confidence.
}
The goal of this experiment was to contrast the
extremal behavior of $\hat f\lc$ vs. $\hat f\hc$.
To this end, we 
numerically optimized
the $\vec p\in[0,1]^n$ so as to maximize
the {\em absolute gap}
\beq
\Delta_n(\vec p):= 
\P(f\lc(\vec X)\neq Y)
-
\P(f\opt(\vec X)\neq Y)
,
\eeq
where 
$f\lc(\vec x)=\sgn\paren{\sum_{i=1}^n (p_i-\hf)x_i}$.
We were surprised to discover that,
though the ratio $
\P(f\lc(\vec X)\neq Y)
/
\P(f\opt(\vec X)\neq Y)
$
can be made arbitrarily large by setting $p_1\approx1$ and the remaining
$p_i<1-\eps$, 
the absolute gap 
appears to be rather small:
we conjecture (with some heuristic justification) that
$\sup_{n\ge1}\sup_{\vec p\in[0,1]^n} \Delta_n(\vec p)=1/16$.
For $\hat f\ba$, we used $\alpha_i=\beta_i=1$ for all $i$. 
The results are reported in Figure~\ref{fig:special5}.
\paragraph{Bayesian setting.
}
In each trial, a vector of expert competences $\vec p\in[0,1]^n$
was drawn independently componentwise, with $p_i\sim\Bdist(1,1)$.
These values (i.e., $\alpha_i=\beta_i\equiv1$) were used for 
$\hat f\ba$.
The results are reported in Figure~\ref{fig:bayes}.
\section{Discussion}
The classic and seemingly well-understood problem of 
the consistency of weighted majority votes continues to reveal
untapped
depth and suggest challenging unresolved questions.
We hope that the results and open problems presented here will 
stimulate future research.
\hide{}
\bibliographystyle{icml2014}
\bibliography{../../mybib}
\appendix
\section*{Appendix: Deferred proofs}
\hide{}
\bepf[Proof of Lemma~\ref{lem:pqw}]
Since $F$ is symmetric about $x=\hf$, it suffices
to prove the claim for $\hf\le x<1$.
We will show that $F$ is concave by examining
its second derivative:
\beq
F''(x) = -\frac{ 2x-1 -2x(1-x)\log(x/(1-x)) }
{ x(1-x)(2x-1)^3}.
\eeq
The denominator is obviously nonnegative on $[\hf,1]$, while the numerator
has the Taylor expansion
\beq
\hide{}
\sum_{n=1}^\infty \frac{2^{2(n+1)}(x-\hf)^{2n+1}}{4n^2-1}\ge0,
\qquad \hf\le x<1
\eeq
(verified through tedious but straightforward calculus).
Since $F$ is concave and symmetric about $\hf$,
its maximum occurs at $F(\hf)=\hf$.
\enpf
\hide{}
\end{document}